\documentclass[letter,10pt]{article}

\usepackage[top=3cm, bottom=3cm, left=3cm, right=3cm, columnsep=0.75cm]{geometry}
\usepackage{caption}
\usepackage{varwidth}
\usepackage{color}
\usepackage{array}
\usepackage{longtable}
\usepackage{algpseudocode}
\algdef{SE}{Begin}{End}{\textbf{begin}}{\textbf{end}}
\usepackage{algorithmicx}
\usepackage{algorithm}
\usepackage{multirow}
\usepackage{graphicx,subcaption}
\usepackage[flushleft]{threeparttable} 
\usepackage{booktabs,caption}
\usepackage{cite}

\usepackage{bmpsize}
\usepackage{multirow}
\usepackage{epsfig}
\usepackage{amssymb}
\usepackage{amsmath}
\usepackage{amsfonts}
\usepackage{amsthm}
\usepackage{mathtools}
\usepackage{multimedia}
\usepackage{tikz}
\usepackage{lipsum}
\usetikzlibrary{positioning}
\usepackage{array,booktabs}
\usepackage{bm}
\usepackage{mathtools}
\usepackage{empheq}
\usepackage{pgfplots}
\usepackage[parfill]{parskip}
\usepackage{txfonts}

\usepackage{fancybox}
\usepackage{anyfontsize} 

\usepackage{sectsty}
\usepackage{subcaption}
\usepackage[utf8]{inputenc}
\usepackage[thinlines,thiklines]{easybmat}
\usepackage{hhline}
\usepackage{arydshln}
\pgfplotsset{compat=1.8}

\makeatletter
\DeclareMathSizes{\f@size}{11}{6}{6}
\makeatother

\DeclareMathOperator{\R}{\varmathbb{R}}
\DeclareMathOperator{\X}{\mathcal{X}}

\DeclareMathOperator{\Z}{\varmathbb{Z}}
\DeclareMathOperator{\E}{\mathbb{E}}

\DeclareMathOperator{\e}{\mathbf e}

\DeclareMathOperator{\x}{\mathbf x}

\DeclareMathOperator{\y}{\mathbf y}

\DeclareMathOperator{\z}{\mathbf z}
\DeclareMathOperator{\w}{\mathbf w}
\DeclareMathOperator{\cc}{\mathbf c}
\DeclareMathOperator{\bb}{\mathbf b}

\DeclareMathOperator{\bz}{\mathbf 0}

\DeclareMathOperator{\eig}{\hbox{eig}}

\newcolumntype{a}{>{\columncolor{BlueGreen}}c}
\newcolumntype{b}{>{\columncolor{Dandelion}}c}
\newcolumntype{d}{>{\columncolor{GreenYellow}}c}
\def\VRHDW#1#2#3{\vrule height #1 depth #2 width #3}
\newcommand{\up}{\VRHDW{1.5em}{0em}{0em}}
\newcommand{\uph}{\VRHDW{1em}{0em}{0em}}

\def\sgn{\mbox{Sign}}

\newcommand{\myeqdl}[1]{\begin{equation} {\begin{array}{rccl} #1 \end{array}}\end{equation}}

\newcommand{\myeq}[1]{$  {#1} $}

\newcommand{\myeqln}[1]{\begin{equation*}  {#1} \end{equation*}}
\newcommand{\myalg}[2]{\begin{equation} {#1 \left\{\begin{array}{l} #2 \end{array}\right.}\end{equation}}

\newcommand{\mycasesn}[2]{\begin{equation*} {#1 \begin{cases}\begin{array}{ll} #2 \end{array}\end{cases}}\end{equation*}}

\theoremstyle{definition}
\newtheorem{theorem}{Theorem}[section]
\newtheorem{corollary}{Corollary}[theorem] 
\newtheorem{lemma}[theorem]{Lemma} 
\newtheorem{assumption}[theorem]{Assumption}
\newtheorem{proposition}[theorem]{Proposition}

\newcommand{\myasmp}[1]{\begin{assumption} {\it #1}\end{assumption}}

\newcommand{\myth}[1]{\begin{theorem} {\it #1}\end{theorem}}
\newcommand{\mylemma}[1]{\begin{lemma} {\it #1}\end{lemma}}
\newcommand{\myprop}[1]{\begin{proposition} {\it #1}\end{proposition}}

\title{On a Randomized Multi-Block ADMM for Solving Selected Machine Learning Problems}
\author{%
 Mingxi Zhu\footnote{Mingxi Zhu is with the Graduate School of Business, Stanford University, USA. Email: mingxiz@stanford.edu}
\and%
Kre\v{s}imir Mihi\'{c}\footnote{Kresimir Mihic is with the School of Mathematics, The University of Edinburgh, UK; and Oracle Labs, Redwood Shores, CA, USA. Email: K.Mihic@sms.ed.ac.uk, kresimir.mihic@oracle.com}
\and Yinyu Ye\footnote{Yinyu Ye is with the Department of Management Science and Engineering, School of Engineering, Stanford University, USA. Email: yyye@stanford.edu.}%
}
\begin{document}

\maketitle

\begin{abstract}

The Alternating Direction Method of Multipliers (ADMM) has now days gained a substantial attention for solving large-scale machine learning and signal processing problems due to the relative simplicity. However, the two-block structure of the classical ADMM still limits the size of the real problems being solved. When one forces a more-than-two-block structure, the convergence speed slows down greatly as observed in practice. Recently, a randomly assembled cyclic multi-block ADMM (RAC-MBADMM) was developed by  the authors for solving general convex quadratic optimization problems where the number of blocks can go greater than two so that each sub-problem has a much smaller size and can be solved much more efficiently. In this paper, we apply this method to solving few selected machine learning and statistic problems related to convex quadratic optimization, such as Linear Regression, LASSO, Elastic-Net, and SVM. We use our  solver to conduct multiple numerical tests, solving both synthetic and large-scale bench-mark problems. Our results show that RAC-MBADMM could significantly outperform other optimization algorithms/codes designed to solve these machine learning problems in both solution time and quality, and match up the performance of the best tailored methods such as Glmnet or LIBSVM. In certain problem regions, e.g., for problems with high dimensional features, RAC-MBADMM also achieves a better performance than that of tailored methods.
\end{abstract}

 \section{Introduction}

We consider the following general convex optimization problems
\myeqdl{
\label{general_model}
\begin{array}{cl}
\min_{\x\in\X}& \sum_{i=1}^{p}f_i(\x_i) + \frac{1}{2} \x^T H \x + \cc^T\x\\
\mbox{s.t.}& \sum_{i=1}^{p}A_i\x_i =\bb
\end{array}
}
where $f_i:\R^{d_i}\mapsto (-\infty,+\infty]$ are closed proper convex functions, $H\in\R^{n\times n}$ is a symmetric positive semi-definite matrix, vector $\cc\in\R^n$. And the problem parameters are the matrix $A=[A_1,\dots,A_p]$, $A_i\in\R^{m\times d_i}$, $i = 1,2,\dots, p$ with $\sum_{i=1}^{p} d_i = n$  and the vector $\bb\in\R^m$. The constraint set $\mathcal X$ is the Cartesian product of possibly non-convex real, closed, nonempty sets, ${\mathcal X} = {\mathcal X_1} \times \dots \times{\mathcal X_p}$, where ${\x_i\in\mathcal X_i} \subseteq \R^{d_i}$. Such problem naturally arises in statistical and learning problems, including elastic net estimation, and supporting vector machine (SVM) problems.

The augmented lagrangian of \ref{general_model} is given by
\myeqdl{\label{eq:lagrangian}\begin{array}{ll}
L_{\mathcal{A}} =&  \sum_{i=1}^{p}f_i(\x_i) + \frac{1}{2} \x^T H \x + \cc^T\x \\
\\
&-y^T(\sum_{i=1}^{p}A_i\x_i - \bb)+\frac{\lambda}{2}||\sum_{i=1}^{p}A_i\x_i - \bb||^2_2 \end{array}}

RAC-MBADMM an iterative algorithm that embeds a Gaussian-Seidel decomposition into each iteration of the augmented Lagrangian method (ALM) (\cite{hestenes:1969,powell:1978}). 
It can be viewed as a decomposition-coordination procedure that decomposes the problem in a random fashion and combines the solutions to small
local sub-problems to find the solution to the original large-scale problem.
The algorithm consists of a cyclic update of randomly constructed blocks of primal variables, $\x_i\in\X_{i}$, followed by a dual ascent type update for Lagrange multipliers $\y$:
\myalg{\label{RAC-ADDM}}{\mbox{Randomly (without replacement) assemble primal}\\
\mbox{variables in $\x$ $^\dagger$ into $p$ blocks $\x_i$, $i=1,\dots,p$,} \\
\mbox{then solve}:\\[0.1cm]
\x_1^{k+1}=\arg\min \{L_{\mathcal{A}}(\x_1,\x_2^k,\dots,\x_p^k,\y^k)\,|\, \x_1\in X_1\}, \\[0.1cm]\vspace{5pt}
\hspace{50pt}\hdots\\
\x_p^{k+1}=\arg\min \{L_{\mathcal{A}} (\x_1^{k+1},\x_2^{k+1},\dots,\x_p\y^k)\,|\, \x_p\in X_p\},\\[0.1cm]
\y^{k+1}=\y^k -\beta(A\x^{k+1}- \bb)
}

RAC-MBADMM can be seen as a generalization of cyclic-ADMM, i.e. cyclic multi-block ADMM is a special case of RAC-MBADMM in which the blocks are constructed at each iteration using a deterministic rule and optimized following a fixed block order. Note that the algorithmic scheme (\ref{RAC-ADDM}) under a deterministic block assembly rule reduces to the classical 2-block ADMM (\cite{glowinski:2014, gabay:1976}) when there are only two blocks ($p=2$) and the
coupled objective vanishes ($H = 0$). 

The classical 2-block ADMM and its convergence have been extensively studied in the literature (e.g. \cite{gabay:1976, eckstein:1992, he:2012, monteiro:2013, deng:2016}). However, the two-block variable structure of the ADMM still limits the practical computational efficiency of the method, because, one factorization of a large matrix is needed at least once even for linear and convex quadratic programming (e.g.,\cite{stellato:2018, zhang:2018}). This drawback may be overcome by enforcing a multi-block structure of the decision variables in the original optimization problem. Due to the simplicity of the latter, there is 
 an active research going on in developing ADMM variants with provable convergence and competitive numerical efficiency and iteration simplicity (e.g. \cite{chen:2017, he:2012, hong:2017, peng:2012}), and on proving global convergence under some special conditions (e.g. \cite{lin:2015,lin:2016,esser:2010,cai:2014}). One of the effective way to resolve divergence issue of multi-block ADMM is to randomly permute the update order block-wise \cite{sunYe:2015}. We denote this algorithm as RP-MBADMM. RP-MBADMM can be seen as a special case of RAC-MBADMM, in which blocks are constructed using some predetermined rule and kept fixed at each iteration, but sub-problems (i.e. blocks minimizing primal variables) are solved in a random order. 

Distributed variants of multi-block ADMM were suggested in \cite{bert:2015, parikh:2014}. The methods convert the multi-block problem into an equivalent two-block problem via variable splitting \cite{bert:1989} and perform a separate augmented Lagrangian minimization over $\x_i$. 
Because of the variable splitting, the distributed ADMM approach increases the number of variables and constraints in the problem, which in turn makes the algorithms not very efficient for large $p$. In addition, the method is not provably working for solving problems with non-separable objective functions.

Our recent paper (\cite{mihic:2019_full}) proved the expected convergence of RAC-MBADMM algorithm. This paper differs from \cite{mihic:2019_full} in the following key aspects
\begin{itemize}
    \item The algorithm implementation in \cite{mihic:2019_full} considers quadratic programming problem with linear constraints, in this paper, we augmented the algorithm for elastic net problem.
    \item The numerical results in \cite{mihic:2019_full} focus mainly on operations problems, including quadratic assignment models, and provide benchmark results of various quadratic programming problems. The focus of this paper is on the statistical and learning models.
\end{itemize}

 We still include some theoretical contributions in \cite{mihic:2019_full} for self-completeness.

The current paper is organized as follows. In the next section we outline
theoretical results with respect to convergence of RAC-MBADMM, followed by 
numerical tests presented in Section \ref{sect:app}. In section \ref{sect:lasso} we compare our general purpose convex quadratic optimization solver \cite{RACQP:code} with glmnet \cite{friedman:2010, friedman:2011}, OSQP \cite{stellato:2018} and Matlab on Linear Regression problems, and in Section \ref{sect:svm} with LIBSVM \cite{Chang:2011} and Matlab on SVM problems. The summary of our contributions with concluding remarks is given
in Section \ref{sect:summary}.

 \begin{flushleft}
\section{Convergence of RAC-MBADMM}
\end{flushleft}
\label{sect:theory}

As shown in \cite{chen:2016}, the convergence result for 2-block ADMM cannot be directly extended to the multi-block case. To remove the possibility of divergence,  \cite{sunYe:2015} shows that with randomly permuted block-wise update order, RP-MBADMM converges in expectation.  

In \cite{chen:2017, chen:2015} the authors focused on solving the linearly constrained convex optimization with coupled convex quadratic objective, and proved the convergence in expectation of RP-MBADMM for the non separable multi-block convex quadratic programming. Proof of convergence in expectation of RAC-MBADMM extends from these results. The detailed proof of the Theorem \ref{th:rac-convergence} is given in \cite{mihic:2019_full} (Section 2.2.2), but outlined in this section for the completeness.

Consider the following linear constrained 
quadratic problem 
\myeqdl{\label{lcqp}
\begin{array}{cl}
     \min & \frac{1}{2}\x^T H \x\ +\ c^T\x  \\[5pt]
     \mbox{s.t.} & A\x \ = \ \bb  \\
\end{array}
}

\vspace{.3cm}
\myasmp{
\label{assumption1}
Assume that for any block of primal variables $\x_i$, 
\myeqln{H_{\sigma_i,\sigma_i} + \beta A^T_{\sigma_i}A_{\sigma_j}\succ0}
where $\sigma_i$ is the index vector describing indices of primal variables of the block $i$.
}
\vspace{.1cm}
\myth{
\label{th:rac-convergence}
Suppose that Assumption $(\ref{assumption1})$ holds, and that RAC-MPADMM (\ref{RAC-ADDM}) is employed to solve problem (\ref{lcqp}). Then the expected output converges to some KKT point of (\ref{lcqp}).
}

Let $\Gamma_{(n,p)}$ denote all possible updating combinations for RAC with $n$ variables and $p$ blocks, and let $\sigma\in\Gamma_{(n,p)}$ denote one specific updating combination for RAC-MBADMM. Then the total  number of updating combinations for RAC-MBADMM is given by
\myeq{|\Gamma_{(n,p)}|=n!/(s!)^{p}}
where $s\in\Z_+$ denotes size of each block with $p\cdot s=n$. 
Let $\upsilon_i\in \Upsilon_{(n,p)}$ denote one specific block composition or partition of $n$ decision variables into $p$ blocks, where $\Upsilon_{(n,p)}$ is the set of all possible block compositions. Then, the total number of all possible block compositions is given by \myeq{|\Upsilon_{(n,p)}|=n!/p!(s!)^{p}}.

Recall the augmented Lagrangian function described with Eq.\ref{eq:lagrangian}, and consider one specific update order generated by RAC-MBADMM, $\sigma\in\Gamma_{RAC(n,p)}$,  $\sigma = [\sigma_1,\dots,\sigma_p]$, where $\sigma_i$ is an index vector of size $s$.
Let $L_{\sigma}\in\R^{n\times n}$ be $s\times s$ block matrix defined with respect to $\sigma_i$ rows and $\sigma_j$ columns as
\mycasesn{(L_{\sigma})_{\sigma_i,\sigma_j}:=}
{H_{\sigma_i,\sigma_j} + \beta A^T_{\sigma_i}A_{\sigma_j},  & \quad i\geq j\\
0,  & \quad \textup{otherwise}
}
and let $R _{\sigma}$ be defined as
\myeqln{R_{\sigma}:= L_\sigma - (H + \beta A^TA).} 
By setting $\z:=(\x;\y)$, RAC-MBADMM could be viewed as a linear system mapping iteration
\myeqln{\z^{k+1}:= M_{\sigma} \z^{k}+\bar{L}^{-1}_{\sigma}\bar{\bb}}
where
\myeqln{
\begin{array}{cccc}
M_{\sigma}:= \bar{L}^{-1}_{\sigma}\bar{R}_{\sigma}, &
\bar{L}_{\sigma}:=\begin{bmatrix}L_{\sigma} & 0\\
\beta A & I\end{bmatrix}, \\
\\
\bar{R}_{\sigma}:=\begin{bmatrix}R_{\sigma} & A^T\\
0 & I\end{bmatrix}, &
\bar{\bb}:=\begin{bmatrix}-\cc+\beta A^T\bb \\
\beta \bb\end{bmatrix} 
\end{array}
}

Finally, the expected mapping matrix $M$ is given by
\myeqln{
M:=\E_{\sigma}(M_\sigma) = \begin{bmatrix}
I-QS & QA^T\\
-\beta A+\beta AQS & I-\beta AQA^T
\end{bmatrix}}
where $S=H + \beta A^TA$ and \myeq{Q:=\E_{\sigma}(L^{-1}_{\sigma})}.

With the preliminaries defined, to prove Theorem \ref{th:rac-convergence}, we follow the same proof structure as in \cite{sunYe:2015, chen:2017, sunYe:2019}, and show that under Assumption \ref{assumption1}:
\begin{enumerate}
 \item[(1)] $\eig(QS)\in[0,\frac{4}{3})$ ;
 \item[(2)] $\forall \lambda\in\eig(M), \eig(QS)\in[0,\frac{4}{3})\implies \|\lambda\|<1$ or $\lambda=1$;
 \item[(3)] if $1\in\eig(M)$, then the eigenvalue 1 has a complete set of eigenvectors;
 \item[(4)] Steps (2) and (3) imply the convergence in expectation of the RAC-MBADMM.
\end{enumerate}

In the proof we make use of Theorem 2 from \cite{chen:2017}, which describes RP-MBADMM convergence in expectation under specific conditions put on matrices $H$ and $A$, and Weyl's inequality, which gives the upper bound on maximum eigenvalue  and the lower bound on minimum eigenvalue of a sum of Hermitian matrices. Proofs for items (2) and (3) are identical to proofs given in \cite{chen:2017}, Section 3.2, so here we focus on proving item (1).
The following lemma completes the proof of expected convergence of RAC.
\vspace{.3cm}
\mylemma{
\label{lemma2}
Under Assumption $\ref{assumption1}$, 
\myeqln{eig(QS)\subset[0,\frac{4}{3})} 
}

To prove Lemma \ref{lemma2}, we first show that for any block structure $\upsilon_i$, the following proposition holds.
\myprop{
\label{Proposition 1} 
$Q_{\upsilon_i}S$ is positive definite and symmetric, and 
\myeqln{eig(Q_{\upsilon_i}S)\subseteq [0,\frac{4}{3})}
}

Then noticing that by definition of $Q$, we have 
\myeq{
QS=\frac{1}{\Upsilon(n,p)}\sum_{\upsilon_i}Q_{\upsilon_i}S
} 
where $Q_{\upsilon_i}S$ is positive definite and symmetric. Let $\lambda_1(A)$ denote the maximum eigenvalue of matrix $A$, then as all $Q_{\upsilon_i}S$ are Hermitian matrices, by Weyl's theorem, we have
\myeqln{
\lambda_1(QS)=\lambda_1(\frac{1}{\Upsilon(d,n)} \sum_{i\in \Upsilon}Q_{\upsilon_i}S)\leq \frac{1}{\Upsilon(d,n)}\sum_{i\in \Upsilon}\lambda_1(Q_{\upsilon_i}S)
}
and as $\lambda_1(Q_{\upsilon_i}S)<\frac{4}{3}$ for each $i$, 
\myeqln{\eig(QS)\subseteq [0,\frac{4}{3})
}
which completes the proof of Lemma \ref{lemma2}, and thus establishes that RAC-MBADMM is guaranteed to converge in expectation. 

Note that expected convergence $\not=$ convergence, but is an evidence for convergence for solving most problems, e.g., when iterates are bounded. 
For the strong notion of convergence we use convergence almost surely as an indicator of RAC-MKADMM convergence. Convergence almost surely as a measure for stability has been used in linear control systems for quite some time. 
\vspace{0.4cm}
\myth{\label{almost surely convergence}
Suppose that Assumption \ref{assumption1} holds, and that RAC-MBADMM (\ref{RAC-ADDM}) is employed to solve problem (\ref{lcqp}). Then the output of RAC-MBADMM converges almost surely to some KKT point of (\ref{lcqp}) if
\myeqln{\rho(\E(M_\sigma\otimes M_\sigma))<1}
where $M\otimes M$ is the Kronecker product of $M$ with itself.
}
In the proof (\cite{mihic:2019_full} Section 2.2.4) we make use of Borel-Cantelli’s theorem to show that RAC-MBADMM converges to the solution almost surely, or a.s. in short.

\section{Solving ML problems}
\label{sect:app}
In this section we apply RAC-MBADMM method and RP-MBADMM method to few selected
machine learning (ML) problems related to convex quadratic optimization, namely Linear
Regression (Elastic-Net) and SVM. To solve the problems we use our general-purpose QP solver, RACQP, which implements RAC-MBADMM algorithm in Matlab \cite{matlab}, and is freely available for download at \cite{RACQP:code}. The experimental results reported in this section  were done using 16-core Intel Xeon CPU E5-2650 machine with 96Gb memory running Debian
linux 4.9.168.

In \cite{mihic:2019_full} we provide numerical evidence that RP-MBADMM suffers from slow convergence to a high precision level on L1-norm of equality constraints. However, the benefit of RP is that it could store and pre-factorize sub-block matrices, as block structure is fixed at each iteration for RP-MBADMM. On the other hand, RAC-MBADMM requires reformulation of sub-blocks for each iteration. In many machine learning problems, including regression and SVM, due to the nature of problem, they require less precision level. This gives benefits to both RAC-MBADMM and RP-MBADMM. Under less precision level on constraints, RAC-MBADMM takes few steps until convergence to target level precision, and time spent on formulation on sub-block is less costly. Similarly, for RP-MBADMM, under less precision level, it could converge within few steps, and will not suffer from slow convergence. In fact, for some regression problems, we fix number of iterations. Under this setup, we would expect that RP-MBADMM performs better in time, compared with RAC-MBADMM.

  \subsection{RAC-MBADMM solution for Linear Regression using Elastic Net}
\label{sect:lasso}
  
For elastic net regression, the model is given by

\myeqdl{\label{eq:el_net_original}
\min_{\beta} \ \frac{1}{2n}(y-\mathbf{X}\beta)^{T}(y-\mathbf{X}\beta) + P_{\lambda,\alpha}(\beta)}

with  \myeqdl{P_{\lambda,\alpha}(\beta)=\lambda\{\frac{1-\alpha}{2}\|\beta\|_2+\alpha\|\beta\|_1\}} 

Model \ref{eq:el_net_original} could be reformulate as the following convex linear constrained problem

\myeqdl{\label{eq:el_net_mod}
\begin{array}{cl}
     \min\limits_{\beta} & \frac{1}{2N}(y-X\beta)^{T}(y-X\beta)+P_{\lambda,\alpha}(z) \\
     \mbox{s.t.} & \beta - \z \ = \ \bz 
\end{array}
}
By adjusting $\alpha$ and $\lambda$, one could obtain different models: for ridge regression, $\alpha = 0$, for lasso $\alpha=1$, and for classic linear regression, $\lambda=0$. And the observations \myeq{\mathbf{X}\in\mathbb{R}^{n\times p }}, where $n$ is number of observations and $p$ is number of features. For the problem to be solved by ADMM, we use variable splitting and reformulate the problem as follows
\myeqdl{\label{eq:el_net}
\begin{array}{cl}
     \min\limits_{\beta} & \frac{1}{2N}(y-X\beta)^{T}(y-X\beta)+P_{\lambda,\alpha}(z) \\
     \mbox{s.t.} & \beta - \z \ = \ \bz 
\end{array}
}

Let $c=-\frac{1}{n}X^Ty$, $A=\frac{X}{\sqrt{n}}$, and let $\gamma$ denote the  augmented Lagrangian penalty parameter with respect to constraint $\beta-z$, and $\xi$ be the dual with respect to constraint $\beta-z$. 
The augmented Lagrangian could then be written as 
\myeqln{
\begin{array}{ll}L_{\lambda}=&\frac{1}{2}\beta^T(A^TA+\gamma I)\beta + (c-\xi)^T\beta+(\xi-\gamma \beta)^Tz+\frac{\gamma}{2} z^Tz+P_{\lambda,\alpha}(z)
\end{array}}

We apply RAC-MBADMM algorithm by partitioning $\beta$ into multi-blocks, but solve $z$ as one block. For any given $\beta_{k+1}$, optimizer $z^*_{k+1}$ has the closed form solution. 
\myeqln{z^*_{k+1}(i)(\beta_{k+1}(i), \xi_k(i))=\frac{S(\xi_k(i)-\gamma\beta_{k+1}(i),\lambda
\alpha)}{(1-\alpha)\lambda+\gamma},}
where $\xi_i$ is the dual variable with respect to constraint $\beta_i-z_i=0$, and $S(a,b)$ is soft-threshold operation \cite{friedman:2007}.
\myeqln{S(a,b)=\begin{cases}
-(a-b),  &\textup{if} \  b<|a|, \ a>0 \\
-(a+b),  &\textup{if} \ b<|a|, \ a\leq 0 \\
0,     & \textup{if}\ b\geq |a|  \\
\end{cases}}

In order to solve classic linear regression directly,  $\mathbf{X}^T\mathbf{X}$ must be positive definite which can not be satisfied for $p>n$. However, RAC-MBADMM only requires  each sub-block  $\mathbf{X}_{sub}^T\mathbf{X}_{sub}$ to be positive definite, so, as long as block size $s<n$, RAC-MBADMM can be used to solve the classic linear regression.

In most experiments, for stopping criteria, we fix the same number of iterations across different solvers. The reason for choosing fixing number of iterations as stopping criteria is mainly because we are focusing on high dimensional data, and setting a high-precision residual tolerance $=||\beta-\z||_1$ as stopping criteria may be too time-consuming. Besides, third-party solvers (e.g. Glmnet) also chooses fixed number of iterations as stopping criteria, provided the tolerance criteria is not meet. We also tested our algorithm under a  high-precision ($10^{-7}$) residual tolerance stopping criteria for relatively medium-sized problems. Our algorithm achieves good performance under both stopping criteria.

\subsubsection{Experimental results}

\underline{Synthetic Data}

In this section we compare our solver with glmnet \cite{friedman:2010, friedman:2011}  and Matlab lasso implementation. 

\underline{Comparison with consensus ADMM}

We first compare our algorithms with consensus ADMM on medium-sized synthetic sparse problems. For consensus ADMM, it applies variable splitting \cite{bert:1989} and introduce auxiliary variables with respect to each block. Although Consensus ADMM possesses beautiful $O(1/k)$ convergence for separable objective functions, these methods do not provably converge for problems with non-separable objectives (e.g. SVM with Gaussian kernel) while ours provably converges linearly in expectation. In this section, we test the performance of consensus ADMM. The problem generated with size of  $X\in\mathbb{R}^{1000\times 100,000}$, $y\in\mathbb{R}^{1000}$ and sparsity of $x$ is $0.95$. 
And here, in stead of fixing number iterations, we use feasibility residual tolerance (set to $10^{-7}$) to be the stopping criteria, and report report run time, number of iterations and the prediction quality $=||X\beta_{alg}-y||_2$. As Glmnet uses coordinate gradient descent without introducing auxiliaries, we cannot fix the same stopping criteria for Glmnet, so here we only compare RAC-MBADMM, RP-MBADMM with consensus ADMM.

\begin{table}[h!]
  \caption{Comparison of RAC-MBBADMM and Consensus-ADMM}
  \label{tbl:Parallel}
  \centering
  \footnotesize
  \begin{tabular}{ccccc}
    \toprule
  &  Num Iter & Time &  Tolerance & Quality\\
  \midrule
RAC  & 31   &56.52&  8.34e-08&1.36e-06\\
RP  & 33   &16.25&  5.15e-08& 6.24e-05\\
Consensus& 470&147.82&9.79e-08&8.48e-07\\
    \bottomrule
  \end{tabular}
\end{table}

From this set of numerical experiments, firstly we found that consensus ADMM performs slower compared with RP-MBADMM or RAC-MBADMM. We observe the similar phenomenon in other numerical tests. Since the consensus ADMM is typically much slower than RAC-MBADMM or RP-MBADMM, we choose not to include its results in the following computational comparisons. Secondly, we observe that although RP takes longer number of iterations to converge, it takes less time, as RP-MBADMM could take advantage on pre-factorization. However, if we fix relative small number of iterations, such advantage is not significant. In fact, in the following numerical experiments, we found that RAC-MBADMM and RP-MBADMM have similar performance on total time.
We compare our solver with glmnet \cite{friedman:2010, friedman:2011}  and Matlab lasso implementation on synthetic data (sparse and dense problems) and benchmark regression data from LIBSVM \cite{Chang:2011}.

\underline{ Synthetic Data}

\underline{ Dense problems}

The data set for dense problems $\mathbf{X}$ is generated uniform randomly with $n=10,000$, $p=50,000$, with zero sparsity, while for the ground truth $\beta^{*}$ we use  standard Gaussian and set sparsity of $\beta^*$ to $0.1$.  Due to the nature of the problem, estimation requires lower feasibility precision, so we fix number of iterations to $10$ and $20$.
Glmnet solver benefits from having a diminishing sequence of $\lambda$, but given that many applications (e.g. see  \cite{mohsen:2015}) require a fixed $\lambda$ value , we decided to use fixed $\lambda$ for all solvers. 
Note that the computation time of RAC-ADMM solver is invariant regardless of whether $\lambda$ is decreasing or fixed.

\begin{table}[h!]
  \centering
  \footnotesize
  \begin{tabular}{ccrrrrrrrr}
    \toprule
\multirow{2}{*}{$\lambda$} & \multirow{2}{*}{Num.} & \multicolumn{4}{c}{Absolute L2 loss} & \multicolumn{4}{c}{Total time [s]} \\ 
\cmidrule(lr){3-6}\cmidrule(lr){7-10}
\uph & iterations & RAC & RP & glmnet & Matlab & RAC & RP & glmnet & Matlab\\ 
\midrule
\multirow{2}{*}{0.01} & 10 & 204.8 & 204.6  & 213.9 & 249.1 & 396.5 &  227.6  & 2465.9 & 1215.2\\ 
 & 20 & 208.1 & 230.2 & 213.9 & 237.1 & 735.2 & 343.9 & 3857.9 & 2218.2\\ 
\uph \multirow{2}{*}{0.1} & 10 & 217.8 &215.6  & 220.5 & 213.1 & 388.7 & 212.5 & 4444.3 &2125.9  \\ 
 & 20 & 272.6  &202.4  &    220.5 & 212.4 & 739.7 & 337.2 &    4452.4 & 2434.6 \\ 

\uph\multirow{2}{*}{1} & 10 & 213.6& 209.0 & 203.1 & 210.5 & 415.3&213.6  & 3021.1 & 1138.9\\ 
 & 20 & 213.8& 212.4 & 210.5 & 203.1 & 686.3& 392.1 & 5295.5 & 1495.6\\
    \bottomrule
  \end{tabular}
  \caption{Comparison on solver performance, dense elastic net model.  Dense problem, $n=10,000,\ p=50,000$}
  \label{tbl:en_synt}
\end{table}

Table \ref{tbl:en_synt} reports on the average cross-validation run-time and the average absolute $L2$ loss for  all possible pairs $(\alpha, \lambda)$ with parameters chosen from  
$\alpha=\{0, \ 0.1, \ 0.2,\dots,1\}$ and $\lambda = \{1,\ 0.01\}$. 
Without specifying, RAC-ADMM solver run-time parameters were identical across the experiments, with augmented Lagrangian penalty parameter $\gamma = 0.1\lambda$ for sparsity $< 0.995$, $\gamma = \lambda$ for sparsity $> 0.995$,  and block size $s=\!=\!100$.

\underline{Large scale sparse problems}

The data set $\mathbf{X}$ is generated uniform randomly with $n=40,000, \ p=4,000,000$.  We set sparsity $=0.998$. For ground truth $\beta^*$ we use standard Gaussian and set sparsity of $\beta^*$ to $0.5$. Similarly, we fix $\lambda$. Notice that from previous experiment, we found that increasing step size from $10$ to $20$ didn't significantly improve prediction error. We further fix number of iteration equals to $10$. 

In Table \ref{tbl:en_synt_1}, we report on the average cross-validation run time and the average absolute $L2$ loss for all possible pairs $(\alpha, \lambda)$ with parameters chosen from  
$\alpha=\{0, \ 0.1, \ 0.2,\dots,1\}$ and $\lambda = \{1,\ 0.01\}$.
We also report best $L2$ loss for each solver.
Notice it take more than $10,000$ s for matlab lasso to solve even one estimation, here we only report comparison between glmnet and RAC. 
\begin{table}[h!]
  \caption{Comparison on solver performance, elastic net model}
  \label{tbl:en_synt_1}
  \centering
  \footnotesize
  \begin{tabular}{ccrrrrrrrrr}
    \toprule
\multirow{2}{*}{$\lambda$} & \multirow{2}{*}{Num.} & \multicolumn{3}{c}{Avg Absolute L2 loss} & \multicolumn{3}{c}{Best Absolute L2 loss}& \multicolumn{3}{c}{Total time [s]} \\ 
\cmidrule(lr){3-5}\cmidrule(lr){6-8} \cmidrule(lr){9-11}
\uph & iterations & RAC & RP & glmnet & RAC &RP & glmnet & RAC & RP & glmnet \\ 
\midrule
\multirow{1}{*}{0.01} & 10 &   1293.3&  1356.7 & 8180.3 & 745.2& 703.52  &   4780.2& 4116.1&2944.8   &17564.2  \\
\uph\multirow{1}{*}{0.1} & 10 &   777.31& 717.92  & 4050.4 &  613.9& 611.79 &3125.6 & 3756.3 &2989.1  & 12953.7\\
\uph\multirow{1}{*}{1} & 10 &   676.17&  671.23  &3124.5 &  615.7&614.79 &1538.9&3697.8&  3003.8 & 8290.5\\ 
    \bottomrule
  \end{tabular}
  \caption*{  Sparse problem, $n=40,000,\ p=4,000,000$}
\end{table}

\underline{Summary of performance}

Experimental results on synthetic data show that RAC-ADMM solver outperforms significantly all other solvers in total time  while being competitive in absolute $L2$ loss. Further RAC-ADMM speedups could be accomplished by fixing block-structure (RP-ADMM).
In terms of run-time, for dense problem, RAC-ADMM is 3 times faster compared with Matlab lasso and 7 times faster compared with glmnet. RP-ADMM is 6 times faster compared with Matlab lasso, and 14 times faster compared with glmnet. For sparse problem, RAC-ADMM is more than 30 times faster compared with Matlab lasso, and 3 times faster compared with glmnet. RP-ADMM is 4 times faster compared with glmnet.

\cite{mihic:2019_full} provided both theoretical and experimental evidence showing that RP-ADMM could be significantly slower compared with RAC-ADMM.RP-ADMM also suffers from slow convergence to a high precision level on L1-norm of equality constraints. However, the benefit of RP-ADMM is that it could store pre-factorized sub-block matrices, as block structure is fixed at each iteration,in contrast to RAC-ADMM which requires  reformulation of sub-blocks at each iteration, what it turn makes each iteration more time-wise costly . In many machine learning problems, including regression, due to the nature of problem, a less precision level is required. This makes RP-ADMM an attractive approach, as it could converge within fewer steps and potentially be faster than RAC-ADMM. 
In addition, while performing simulations we observed that increasing number of iteration does not significantly improve performance of prediction. In fact, absolute $L2$ loss remains similar even when number of iteration is increased to 100. This further gives an advantage to RP-ADMM, as it benefits the most when number of iteration is relatively small. 

\underline{Benchmark instances}

In section we report on experiments on E2006-tfidf and log1pE2006 regression data from LIBSVM \cite{Chang:2011}. 

\underline{E2006-tfidf}

For E2006-tfidf, feature size is $150,360$, number of training data is $16,087$ and number of testing data is $3,308$. The null training error of test set is $221.8758$. Note that the sparsity of this data set is $0.991$, and our RAC-MBADMM have competitive performance in both dense synthetic problem and sparse problem.
In this setup, we fix number of iterations equals to $10$. Similarly, we vary $\lambda = \{1, \ 0.01\}$, and $\alpha = \{0, \ 0.1, \ 0.2,\dots,1\}$. We use the training set to predict $\beta^*$, and compare the model error (ME) of test set across different solvers. 

Note that for two-block ADMM based solver, including OSQP \cite{stellato:2018} and Matlab lasso, due to the inefficiency of factorizing a big matrix, for a single value of $\alpha$, both solver takes more than $1000$ seconds to solve the problem for even $10$ iterations. On the other hand, glmnet, which is similar to our algorithm using a cyclic coordinate descent algorithm on each variables, performs significantly faster than OSQP and matlab lasso. However, glmnet can still be inefficient, as a complete cycle through all $p$ variables requires $O(pN)$ operations \cite{friedman:2010}. However, for RAC-MBADMM, we only need to calculate $\mathbf{X}^T_{sub}\mathbf{X}_{sub}$ online, whenever a sub block is selected, and solve the small $s\times s$ sized problem each time.  With RAC-MBADMM, we do not need to factorize a huge matrix, nor did we need to calculate a big matrix $\mathbf{X}^T\mathbf{X}$, which could be time and memory-consuming. RAC-MBADMM requires less memory and is more computationally efficient.

In Table \ref{tbl:en_ben_l} we show the performance of OSQP and matlab lasso for $\alpha=1$, and in Table \ref{tbl:en_ben} we compare our solver with glmnet, by averaging over the model error of test set for different $\alpha$. The table shows the average runtime and training error collected from experiments with $\alpha=\{0,0.1,\dots,1\}$. Here we pre-factorize matrix for RP.

\begin{table}[h!]
\caption{E2006-tfidf Lasso Problem Performance Summary}
\label{tbl:en_ben_l}
  \centering
  \footnotesize
  \begin{tabular}{ccc}
    \toprule
  Solver &   Training ME &Total time [s]\\ 
\midrule
OSQP& 64.0 & 1482.5 \\
matlab  & 61.1 & 3946.6\\
\bottomrule
\end{tabular}
\caption*{\footnotesize{$\alpha =1$, $\lambda = 0.01$}}
\vspace{-10pt}
\end{table}

\begin{table}[h!]
  \caption{E2006-tfidf performance summary }
  \label{tbl:en_ben}
  \centering
  \footnotesize
  \begin{tabular}{ccccccc}
    \toprule
\multirow{2}{*}{$\lambda$}   &  \multicolumn{3}{c}{Training ME} & \multicolumn{3}{c}{ Total time [s]}  \\ 
\cmidrule(r){2-4}\cmidrule(r){5-7}
 \uph &  RAC  & RP &  glmnet  &  RAC &RP &  glmnet\\ 
\midrule
0.01  &  22.4& 22.4 &  29.9  &  106.5& 50.9 &  653.2\\ 
0.1  &  22.1&  22.1&  22.7  &   100.5& 51.9  &  269.3\\
1  &  25.7&25.7  &  23.5  &  102.5  &54.2 & 282.9 \\
    \bottomrule
  \end{tabular}
\vspace{-10pt}
\end{table}

We observe that RAC-MBADMM solver is faster compared to glmnet for all different parameters and that it achieves the best training model error, $22.0954$, among all the solvers. 

In terms of time, RAC-MBADMM is 14 times faster compared with OSQP, 38 times faster compared with matlab lasso, and 4 times faster compared with glmnet. RP-MBADMM is 28 times faster compared with OSQP, 18 times faster compared with matlab lasso and 8 times faster compared with glmnet.

\underline{log1pE2006}

For log1pE2006, feature size is $4,272,227$, number of training data is $16,087$ and number of testing data is $3,308$. The null training error of test set is $221.8758$. The sparsity of this data set is $0.998$.

In Table \ref{tbl:en_ben_log_1} we show the performance of OSQP and matlab lasso for $\alpha=1$, and in Table \ref{tbl:en_ben_log} we compare our solver with glmnet, by averaging over the model error of test set for different $\alpha$. The table shows the average runtime and training error collected from experiments with $\alpha=\{0,0.1,\dots,1\}$. For this set of experiment, we didn't pre-factorize for RP as pre-factorization on matrix is more time consuming.

\begin{table}[h!]
\caption{log1pE2006 Lasso Problem Performance Summary}
\label{tbl:en_ben_log_1}
  \centering
  \footnotesize
  \begin{tabular}{ccc}
    \toprule
  Solver &   Training ME &Total time [s]\\ 
\midrule
OSQP&  66.6  &  11437.4\\
matlab  & - & more than 3 days\\
\bottomrule
\end{tabular}
\caption*{\footnotesize{$\alpha =1$, $\lambda = 0.01$}}
\vspace{-10pt}
\end{table}

\begin{table}[h!]
  \caption{log1pE2006 performance summary }
  \label{tbl:en_ben_log}
  \centering
  \footnotesize
  \begin{tabular}{ccccccc}
    \toprule
\multirow{2}{*}{$\lambda$}   &  \multicolumn{3}{c}{Training ME} & \multicolumn{3}{c}{ Total time [s]}  \\ 
\cmidrule(r){2-4}\cmidrule(r){5-7}
 \uph &  RAC& RP  &  glmnet  &  RAC& RP  &  glmnet\\ 
\midrule
0.01  & 43.0 &  41.8& 22.0  &  962.2&   722.5 &  7639.6\\ 
0.1  &  30.8&   31.8&  22.5  &  978.7&  721.4 &  4945.2\\ 
1  &    32.1& 35.5 &  29.3  &  958.5&  749.2  &  1889.5\\
    \bottomrule
  \end{tabular}
\vspace{-10pt}
\end{table}

In terms of model error, RAC-MBADMM and RP-MBADMM are still competitive and is of same level comapared with glmnet. Both glmnet and RAC-MBADMM outperforms OSQP in terms of model error.

The results show that RAC-ADMM and RP-ADMM are still competitive and are of same level as glmnet with respect to model error, and all outperform OSQP and Matlab. In terms of run-time, RAC-ADMM is 12 times faster than OSQP, and 5 times faster than glmnet. RP-ADMM is 16 and 7 times faster than OSQP and glmnet, respectively.

\subsection{RAC-MBADMM solution for SVM}
\label{sect:svm}

A Support Vector Machine (SVM) is a machine learning method for classification, regression, and other learning tasks. The method learns a mapping between the features $\x_i\in\R^r$, $i=1,\dots n$ and the target label $y_i\in\{-1,1\}$ of a set of data points using a {\it training set} and constructs a hyperplane \myeq{\w^T\phi(\x)+b} that separates the data set. 
This hyperplane is then used to predict the class of further data points.  
The objective uses Structural Risk Minimization principle which aims to minimize the empirical risk (i.e. misclassification  error) while maximizing the confidence interval (by maximizing the separation margin) \cite{vapnik:1998, vapnik:1995}.

Training an SVM is a convex optimization problem, with multiple formulations, such as C-support vector classification (C-SVC), $\upsilon$-support vector
classification ($\upsilon$-SVC), $\epsilon-$support vector regression ($\epsilon-$SVR), and many more. As our goal is to compare RACQP, a general QP solver, with specialized SVM software and not to compare SVM methods themselves, we decided on using C-SVC (\cite{boser:1992, cortes:1995}), with the dual problem formulated as
\myeqdl{\label{svm}
\begin{array}{cc}
     \min\limits_{\z} & \frac{1}{2}\z^T Q \z\ -\ \e^T\z  \\
     \mbox{s.t.} & \y^T\z \ = \ 0  \\
          & \z\in[0,C]
\end{array}
}
 
with $Q\!\in\!\R^{n\times n}$, $Q\!\succeq\!0$, \myeq{q_{i,j}=y_iy_jK(\x_i,\x_j)}, where \myeq{K(\x_i,\x_j):=\phi(\x_i)^T\phi(\x_j)} is a kernel function, and regularization parameter $C\!>\!0$. The optimal $\w$ satisfies \myeq{\w=\sum_{i=1}^ny_i\z_i\phi(\x_i)},  and 
the bias term $b$ is calculated using the support vectors that lie on the margins (i.e. $0<\z_i<C$) as \myeq{b_i= \w^T\phi(\x_i) - y_i }.  To avoid numerical stability issues, $b$ is then found by averaging over $b_i$.  The decision function is defined with \myeq{f(\x) =\sgn(\w^T\phi(\x)+b)}.

We compare RACQP with LIBSVM \cite{Chang:2011}, due its popularity, and with Matlab-SVM , due to its ease of use. These methods implement specialized approaches to address the SVM problem (e.g. LIBSVM uses a Sequential Minimal Optimization, SMO, type decomposition method  \cite{fan:2005, bottou:2007}), while our approach solves the optimization problem (\ref{svm}) directly. 

The LIBSVM benchmark library provides a large set of instances for SVM, and we selected a representative subset: training data sets with sizes ranging from 20,000 to 580,000; number of features from eight to 1.3 million. We use the test data sets when provided, otherwise, we create test data by randomly choosing 30\% of testing data and report cross-validation accuracy results.

In Table \ref{tbl:svm} we report on model training run-time and accuracy, defined as (num. correctly predicted data)/(total testing data size)$\times$100\%. RAC-ADMM parameters were as follows: max block size $s=100, 500,$ and $1000$ for small, medium and large instances, respectively and augmented Lagrangian penalty $\beta=0.1p$, where $p$ is the number of blocks, which in this case is found to be $p=\lceil n/s \rceil$ with $n$ being the size of training data set. 
In the experiments we use Gaussian  kernel, \myeq{K(\x_i,\x_j)\!=\!\exp(-\frac{1}{2\sigma^2}\|\x_i-\x_j\|^2)}. 
Kernel parameters $\sigma$ and $C$ were estimated by running a {\it grid-check} on  cross-validation. We tried different pairs $(C,\sigma)$ and picked those that returned the best cross-validation accuracy (done using randomly choose 30\% of train data) when instances were solved using RAC-ADMM. Those pairs were then used to solve the instances with LIBSVM and Matlab. The pairs were chosen from a relatively coarse grid, $\sigma,C\in\{0.1, 1, 10\}$ because the goal of this experiment is to compare RAC-ADMM with heuristic implementations rather than to find the best classifier.
Termination criteria were either primal/dual residual tolerance ($\epsilon_p = 10^{-1}$ and $\epsilon_d = 10^{-0}$) or maximum number of iterations, $k=10$, whichever occurs the first. Dual residual was set to such a low value because empirical observations showed that restricting the dual residual does not significantly increase accuracy of the classification but effects run-time disproportionately. Maximum run-time was limited to 10 hours for mid-size problems, and unlimited for the large ones. Run-time shown in seconds, unless noted otherwise.

\begin{table}[h!]
  \centering
  \footnotesize
\begin{threeparttable}
  \begin{tabular}{lllrrrrrrr}
    \toprule
 \multirow{2}{*}{Instance} & \multirow{2}{*}{Training} 
&\multirow{2}{*}{Testing}   
& \multirow{2}{*}{Num.} & \multicolumn{3}{c}{Accuracy [\%]}  &\multicolumn{3}{c}{Training run-time [s]}  \\ 
\cmidrule(r){5-7}\cmidrule(r){8-10}
\uph name &  set size  
&  set size 
&  features  &  RAC  &  LIBSVM  &  Matlab  &  RAC  &  LIBSVM  &  Matlab\\ 
\midrule
\up
a8a & 22696 & 9865 & 122 & 76.3 & 78.1 & 78.1 & 91 & 250 & 2653\\ 
w7a & 24692 & 25057 & 300 & 97.1 & 97.3 & 97.3 & 83 & 133 & 2155\\ 
rcv1.binary & 20242 & 135480 & 47236 & 73.6 & 52.6 & -- & 78 & 363 & 10+h\\ 
news20.binary$^*$ & 19996 & 5998 & 1355191 & 99.9 & 99.9 & -- & 144 & 3251 & NA\\ 
\up
a9a & 32561 & 16281 & 122 & 76.7 & 78.3 & 78.3 & 211 & 485 & 5502\\ 
w8a & 49749 & 14951 & 300 & 97.2 & 99.5 & 99.5 & 307 & 817 & 20372\\ 
ijcnn1 & 49990 & 91701 & 22 & 91.6 & 91.3 & 91.3 & 505 & 423 & 0\\ 
cod\_rna & 59535 & 271617 & 8 & 79.1 & 73.0 & 73.0 & 381 & 331 & 218\\  
real\_sim$^*$ & 72309 & 21692 & 20958 & 69.5 & 69.5 & -- & 1046 & 9297 & 10+h\\ 
\up
skin\_nonskin$^*$ & 245057 & 73517 & 3 & 99.9 & 99.9 & -- & 2.6h & 0.5h & NA\\  
webspam\_uni$^*$ & 350000 & 105000 & 254 & 64.3 & 99.9 & -- & 13.8h & 11.8h & NA\\ 
covtype.binary$^*$ & 581012 & 174304 & 54 & 91.3 & 99.9 & -- & 16.2h & 45.3h & NA\\
    \bottomrule
  \end{tabular}
\begin{tablenotes}
  \item[*]  No test set provided, using 30\% of randomly chosen data from the training set. Reporting cross-validation accuracy results.
  \end{tablenotes}
\end{threeparttable}
  \caption{Model training performance comparison for SVM }
  \label{tbl:svm}
\end{table}

The results show that RACQP produces classification models of competitive quality as models produced by specialized software implementations in a much shorter time.  RACQP is in general faster than LIBSVM (up to 27x) except for instances where ratio of number of observations $n$ with respect to number of features $r$ is very large. It is noticeable that while producing (almost) identical results as LIBSVM, the Matlab implementation is significantly slower.

For small and mid-size instances (training test size $<$ 100K) we tried, the difference in accuracy prediction is less than 2\%, except for problems where test data sets are much larger than the training sets. In the case of ``rcv1.binary'' instance test data set is 5x larger than the training set, and for ``cod\_rna'' instance is 4x larger. In both cases RACQP outperforms LIBSVM (and Matlab) in accuracy, by 20\%  and  9\%, respectively. 

All instances except for ``news20.binary'' have $n>>r$ and the choice of the Gaussian kernel is the correct one. For instances where the number of features is larger than the number of observations, linear kernel is usually the better choice as the separability of the model can be exploited \cite{woodsend:2011} and problem solved to similar accuracy in a fraction of time required to solve it with the non-linear kernel. 
The reason we used the Gaussian kernel on ``news20.binary' instance is that we wanted to show that RACQP is only mildly affected by the feature set size. Instances of similar sizes but different number of features  are all solved by RACQP in approximately the same time, which is 
in contrast with  LIBSVM and Matlab that are both affected by the feature space size. LIBSVM slows down significantly while Matlab, in addition to slowing down could not solve ''news.binary`` -- the implementation of fitcsvm() function that invokes Matlab-SVM algorithm requires full matrices to be provided as the input which in the case of ''news.binary`` requires 141.3GB of main memory.

``Skin\_nonskin'' benchmark instance ``marks'' a point where our direct approach starts showing weaknesses -- LIBSVM is 5x faster than RACQP because of the fine-tuned heuristics which exploit very small feature space (with respect to number of observations). The largest instance we addressed is ``covtype.binary'', with more than half of million observations and the (relatively) small feature size ($p=54$). For this instance, RACQP continued slowing down proportionately  to the increase in problem size, while LIBSVM experienced a large hit in run-time performance, requiring almost two days to solve the full size problem. This indicates that the algorithms employed by LIBSVM are put to the limit and specialized algorithms (and implementations) are needed to handle large-scale SVM problems. 
RACQP accuracy is lower than that of LIBSVM, but can be improved by tightining residual tolerances under the cost of increased run-time.

For large-size problems RACQP performance degraded, but the success with the mid-size problems suggests that a specialized  ``RAC-SVM'' algorithm could be developed to address very large problems.
Such a solution could merge RAC-ADMM algorithm 
with heuristic techniques to (temporarily) reduce the size of the problem (e.g. \cite{joachims:1998}),  smart 
kernel approximation techniques, probabilistic approach(es) to shrinking the support vector set (e.g. \cite{rudi:2017}), and similar.

 \section{Summary}
\label{sect:summary}

We apply a general purpose convex quadratic optimization solver, randomly assembled cyclic multi-block
ADMM (RAC-MBADMM), to solving few selected machine learning problems such as Linear Regression, LASSO, Elastic-Net, and SVM. Our preliminary numerical tests, solving
both synthetic and large-scale bench-mark problems, indicate that our solver significantly outperforms other optimization algorithms/codes
designed to solve these machine learning problems in both solution time and quality. Our solver also matches the performance of the best tailored methods such as Glmnet and LIBSVM, or often gives better results than that of tailored methods. In addition, our solver uses much less computation memory space than other ADMM based method do, so that it is suitable in real applications with big data.

There are lots of interesting extensions in both theoretical and numerical aspects. Theoretically, for elastic-net type of problems, proving the expected convergence rate when observations $X$ is a random matrix is an important extension. Numerically, there exists clear trade-off on block-size/ block-number selection. Larger blocks implies less number of iterations, however factorization on larger blocks takes more time. More experimentation could be done on selection of optimal block-size. 
\bibliography{references}
\bibliographystyle{siam}  

\end{document}